\documentclass[aoas,preprint]{imsart}

\RequirePackage[OT1]{fontenc}
\RequirePackage{amsthm,amsmath}
\RequirePackage[numbers]{natbib}
\RequirePackage[colorlinks,citecolor=blue,urlcolor=blue]{hyperref}
\usepackage[pdftex]{graphicx}

\startlocaldefs
\numberwithin{equation}{section}
\theoremstyle{plain}

\endlocaldefs

\begin{document}

\begin{frontmatter}
\title{A new approach\\ to Poissonian two-armed bandit problem} \runtitle{Poissonian two-armed bandit}

\begin{aug}
\author{\fnms{Alexander} \snm{Kolnogorov}\thanksref{m1}\ead[label=e1]{Alexander.Kolnogorov@novsu.ru}},

\runauthor{A. Kolnogorov}

\affiliation{Yaroslav-the-Wise Novgorod State
University\thanksmark{m1}}

\address{41 B.Saint-Petersburgskaya Str.,
Velikiy Novgorod, Russia, 173003\\
Applied Mathematics and Information Science Department\\
\printead{e1}}

\end{aug}

\begin{abstract}
We consider a continuous time two-armed bandit problem in which
incomes are described by Poissonian processes. We develop Bayesian
approach with arbitrary prior distribution. We present two
versions of recursive equation for determination of Bayesian
piece-wise constant strategy and Bayesian risk and partial
differential equation in the limiting case. Unlike the previously
considered Bayesian settings our description uses current history
of the process and not evolution of the posterior distribution.
\end{abstract}



\end{frontmatter}

\def \mE{\mathrm{E}}
\def \tR{\tilde{R}}

\section{Introduction}\label{intro}

We consider a continuous time two-armed bandit problem. This
setting results either in Poissonian  or in a diffusion two-armed
bandit. Quite general Poissonian two-armed bandit was considered
in~\cite{PS1,PS2}. In \cite{Mand} consideration of Poissonian and
diffusion bandit problems is restricted to the case of independent
arms and discounted rewards. An interesting though a special case
of diffusion two-armed bandit is presented in \cite{BF}. Some
approaches to a discrete time two-armed bandit problem are
presented in \cite{Sragovich}, \cite{CBL}, \cite{Koln18}. In the
present article, we develop a new general approach to Poissonian
two-armed bandit in Bayesian setting.

Formally, Poissonian two-armed bandit is a continuous-time random
controlled process $X(t)$. Its values are usually interpreted as
incomes  and depend only on chosen actions $y(t)$ as follows. If
on the time interval $t' \in [\tau, \tau+t]$, $t>0$ the action
$y(t')=\ell$ was chosen then
\begin{gather}\label{a1}
\Pr\left(X(\tau+t)-X(\tau\right)=i)=p(i,t;\lambda_\ell)=\frac{(\lambda_\ell
t)^i}{i!} e^{-\lambda_\ell t}, \quad i=0,1,2,\dots
\end{gather}
$\ell=1,2$. Thus a vector parameter $\theta=(\lambda_1,\lambda_2)$
completely describes considered Poissonian two-armed bandit. The
set of admissible values of parameters $\Theta$ is assumed to be
known.

A control strategy generally assigns a random choice of the action
at the point of time $t$ depending on currently observed history
of the process, i.e. cumulative times of both actions applications
$t_1, t_2$ ($t_1+t_2=t$) and corresponding cumulative incomes
$X_1, X_2$. In what follows, current values $X_1, X_2$ at the
point of time $t$ are denoted by $X_1(t), X_2(t)$. If one knew
$\lambda_1,\lambda_2$, he should always choose the action
corresponding to the largest of them, his total expected income on
the control horizon $T$ would thus be equal to
$T\max(\lambda_1,\lambda_2)$. But if he uses some strategy
$\sigma$, his total expected income is less than maximal by the
value
\begin{gather}\label{a2}
L_T(\sigma,\theta)=T
\max(\lambda_1,\lambda_2)-\mE_{\sigma,\theta}(X_1(T)+X_2(T))
\end{gather}
which is called the regret. Here $\mE_{\sigma,\theta}$ denotes the
mathematical expectation with respect to the measure generated by
strategy $\sigma$ and parameter $\theta$.

Let's assign a prior distribution density
$\mu(\theta)=\mu(\lambda_1,\lambda_2)$ on the set of parameters
$\Theta$. Corresponding Bayesian risk is defined as follows
\begin{gather}\label{a3}
R_T(\mu)= \inf_{\{\sigma\}}\int_\Theta
L_T(\sigma,\theta)\mu(\theta) d\theta,
\end{gather}
the optimal strategy $\sigma^B$ is called Bayesian strategy. The
minimax risk on the set $\Theta$ is defined as
\begin{gather}\label{a4}
R^M_T(\Theta)= \inf_{\{\sigma\}}\sup_\Theta L_T(\sigma,\theta),
\end{gather}
corresponding optimal strategy $\sigma^M$ is called minimax
strategy.

A direct method of determining minimax strategy and minimax risk
does not exist. However, one can determine them with the use of
the main theorem of the theory of games. According to this theorem
the following equality holds
\begin{gather}\label{a5}
R^M_T(\Theta)= R_T(\mu_0)= \sup_{\{\mu\}}R_T(\mu),
\end{gather}
i.e. minimax risk is equal to the Bayesian one calculated with
respect to the worst-case prior distribution and minimax strategy
coincides with corresponding Bayesian strategy. Note that in case
of finite set  $\Theta$ determination of the minimax risk
according to equality \eqref{a5} is not laborious because Bayesian
risk is a concave function of the prior distribution.

The rest of the paper is organized as follows. Recursive
Bellman-type equation for determining  Bayesian risk for
piece-wise constant strategies is presented in Section~\ref{equ1}.
Note that our approach differs from presented in~\cite{PS1},
\cite{PS2} because we recalculate Bayesian risk with respect to
current statistics $(X_1,t_1,X_2,t_2)$ and in~\cite{PS1},
\cite{PS2} recalculations are implemented with respect to current
posterior distribution and $t=t_1+t_2$. Our approach is applied to
quite general sets $\Theta$. The approach presented in~\cite{PS1},
\cite{PS2} is applied to finite sets of parameters and
generalization to arbitrary sets is not obvious. In
Section~\ref{equ2}, another version of recursive equation is
derived. In a limiting case, we obtain a partial differential
equation which is presented in Section~\ref{lim}.

\section{Recursive equation}\label{equ1}

Let's consider piece-wise constant strategies
$\{\sigma_\ell(X_1,t_1,X_2,t_2)\}$. To this end, we assume that
control horizon is partitioned into a number of intervals of the
length $\Delta$ on which the chosen action does not change. Hence,
$T=N\Delta$ and for any  $n_1+n_2=n<N$, $t_1=n_1 \Delta$, $t_2=n_2
\Delta$  we have $\Pr(y(t')=\ell)=\sigma_\ell(X_1,t_1,X_2,t_2)$
where $\sigma_\ell(X_1,t_1,X_2,t_2)$ is constant on the time
interval $t' \in [n\Delta, (n+1)\Delta]$. The posterior
distribution at the point of time $t=t_1+t_2$ is calculated as
\begin{gather}\label{a6}
\mu(\lambda_1,\lambda_2|X_1,t_1,X_2,t_2)=\frac{p(X_1,t_1;\lambda_{1})p(X_2,t_2;\lambda_{2})
\mu(\lambda_1,\lambda_2)}{\mu(X_1,t_1,X_2,t_2)},
\end{gather}
where
\begin{gather}\label{a7}
\mu(X_1,t_1,X_2,t_2)=\iint_\Theta
p(X_1,t_1;\lambda_{1})p(X_2,t_2;\lambda_{2})\mu(\lambda_1,\lambda_2)d\lambda_1
d\lambda_2.
\end{gather}

Since $p(0,0;\lambda)=1$, this formula remains correct if $t_1=0$
and/or $t_2=0$. Denote $x^+=\max(x,0)$. With the use of \eqref{a1}
we obtain the following standard recursive Bellman-type equation
for determining Bayesian risk \eqref{a3} with respect to the
posterior distribution \eqref{a6}
\begin{gather}\label{a8}
R(X_1,t_1,X_2,t_2)= \min (R^{(1)}(X_1,t_1,X_2,t_2),
R^{(2)}(X_1,t_1,X_2,t_2)),
\end{gather}
where
\begin{gather}\label{a9}
R^{(1)}(X_1,t_1,X_2,t_2)= R^{(2)}(X_1,t_1,X_2,t_2)=0
\end{gather}
if $t_1+t_2=T$ and then
\begin{gather}\label{a10}
\begin{array}{c}
R^{(1)}(X_1,t_1,X_2,t_2)=\displaystyle{\iint_\Theta}\mu(\lambda_1,\lambda_2|X_1,t_1,X_2,t_2)\\
\times \Big( (\lambda_{2}-\lambda_{1})^+\Delta
+\displaystyle{\sum_{j=0}^\infty} R(X_1+j,t_1+\Delta,X_2,t_2)
p(j,\Delta;\lambda_{1}) \Big)
d\lambda_1 d\lambda_2,\\
R^{(2)}(X_1,t_1,X_2,t_2)=\displaystyle{\iint_\Theta}\mu(\lambda_1,\lambda_2|X_1,t_1,X_2,t_2)\\
 \times \Big( (\lambda_{1}-\lambda_{2})^+\Delta
+\displaystyle{\sum_{j=0}^\infty} R(X_1,t_1,X_2+j,t_2+\Delta)
p(j,\Delta;\lambda_{2}) \Big) d\lambda_1 d\lambda_2.
\end{array}
\end{gather}
Here  $\{R^{(\ell)}(X_1,t_1,X_2,t_2)\}$ are expected losses if
initially the $\ell$-th action is applied at the control horizon
of the length $\Delta$ and then control is optimally implemented
($\ell=1,2$).

Bayesian risk \eqref{a3} is calculated by the formula
\begin{gather}\label{a11}
   R_T(\mu)=R(0,0,0,0).
\end{gather}

Equation \eqref{a8}--\eqref{a10} determine at the same time
Bayesian risk and Bayesian strategy. Bayesian strategy prescribes
to choose $\ell$-th action (i.e $\sigma_\ell(X_1,t_1,X_2,t_2)=1$)
if $R^{(\ell)}(X_1,t_1,X_2,t_2)$ has smaller value. In case of a
draw $R^{(1)}(X_1,t_1,X_2,t_2)=R^{(2)}(X_1,t_1,X_2,t_2)$ the
choice is arbitrary.

\section{Another version of recursive equation}\label{equ2}

In this section, we obtain another version of recursive
Bellman-type equation. Let's denote
\begin{gather*}
   \tR(X_1,t_1,X_2,t_2)=R(X_1,t_1,X_2,t_2)\times
   \mu(X_1,t_1,X_2,t_2),
\end{gather*}
where $\{R(X_1,t_1,X_2,t_2)\}$ are Bayesian risks calculated with
respect to the posterior distribution~\eqref{a6} and
$\{\mu(X_1,t_1,X_2,t_2)\}$ are defined in~\eqref{a7}. Then the
following recursive equation holds
\begin{gather}\label{a12}
\tR(X_1,t_1,X_2,t_2)= \min (\tR^{(1)}(X_1,t_1,X_2,t_2),
\tR^{(2)}(X_1,t_1,X_2,t_2)),
\end{gather}
where
\begin{gather}\label{a13}
\tR^{(1)}(X_1,t_1,X_2,t_2)= \tR^{(2)}(X_1,t_1,X_2,t_2)=0
\end{gather}
if $t_1+t_2=T$ and then
\begin{gather}\label{a14}
\begin{array}{c}
\tR^{(1)}(X_1,t_1,X_2,t_2)=g^{(1)}(X_1,t_1,X_2,t_2)\times \Delta\\
+\displaystyle{\sum_{j=0}^\infty}
\tR(X_1+j,t_1+\Delta,X_2,t_2)\times
\frac{t_1^{X_1}\Delta^j(X_1+j)!}
{(t_1+\Delta)^{X_1+j}X_1! j!} ,\\
\tR^{(2)}(X_1,t_1,X_2,t_2)=g^{(2)}(X_1,t_1,X_2,t_2)\times \Delta\\
+\displaystyle{\sum_{j=0}^\infty} \tR(X_1,t_1,X_2+j,t_2+\Delta)
\times \frac{t_2^{X_2}\Delta^j(X_2+j)!} {(t_2+\Delta)^{X_2+j}X_2!
j!},
\end{array}
\end{gather}
where
\begin{gather*}
g^{(1)}(X_1,t_1,X_2,t_2)=\iint_\Theta(\lambda_{2}-\lambda_{1})^+
p(X_1,t_1;\lambda_{1})p(X_2,t_2;\lambda_{2})\mu(\lambda_1,\lambda_2) d\lambda_1 d\lambda_2,\\
g^{(2)}(X_1,t_1,X_2,t_2)=\iint_\Theta(\lambda_{1}-\lambda_{2})^+
p(X_1,t_1;\lambda_{1})p(X_2,t_2;\lambda_{2})\mu(\lambda_1,\lambda_2)
d\lambda_1 d\lambda_2.
\end{gather*}
Bayesian strategy prescribes to choose  $\ell$-th action (i.e
$\sigma_\ell(X_1,t_1,X_2,t_2)=1$) if
$\tR^{(\ell)}(X_1,t_1,X_2,t_2)$ has smaller value. In case of a
draw $\tR^{(1)}(X_1,t_1,X_2,t_2)=\tR^{(2)}(X_1,t_1,X_2,t_2)$ the
choice is arbitrary. Bayesian risk~\eqref{a3} is calculated by the
formula
\begin{gather}\label{a15}
   R_T(\mu)=\tR(0,0,0,0).
\end{gather}

Formulas~\eqref{a12}--\eqref{a15} follow
from~\eqref{a8}--\eqref{a11}. One should multiply left-hand side
and right-hand side of \eqref{a10} by $\mu(X_1,t_1,X_2,t_2)$ and
implement mathematical transformations.

\section{A limiting description}\label{lim}

In this section, we  consider the case when $\Delta$ has a small
value. In this case \eqref{a14} takes the form
\begin{gather}\label{a16}
\begin{array}{c}
\tR^{(1)}(X_1,t_1,X_2,t_2)=g^{(1)}(X_1,t_1,X_2,t_2) \Delta\\
+ \tR(X_1,t_1+\Delta,X_2,t_2) -\tR(X_1,t_1+\Delta,X_2,t_2)X_1
t_1^{-1}\Delta\\
+\tR(X_1+1,t_1+\Delta,X_2,t_2)(X_1+1) t_1^{-1}\Delta+o(\Delta),\\
\tR^{(2)}(X_1,t_1,X_2,t_2)=g^{(2)}(X_1,t_1,X_2,t_2) \Delta\\
+ \tR(X_1,t_1,X_2,t_2+\Delta) -\tR(X_1,t_1,X_2,t_2+\Delta)X_2
t_2^{-1}\Delta\\
+ \tR(X_1,t_1,X_2+1,t_2+\Delta)(X_2+1)
t_2^{-1}\Delta+o(\Delta),
\end{array}
\end{gather}
Equation~\eqref{a16} must be complemented with~\eqref{a12} which
now is written as
\begin{gather}\label{a17}
\min_{\ell=1,2} \left(
\tR^{(\ell)}(X_1,t_1,X_2,t_2)-\tR(X_1,t_1,X_2,t_2)\right)=0.
\end{gather}
By~\eqref{a16}--\eqref{a17} one derives in the limiting case (as
$\Delta \to +0$) the following partial differential equation
\begin{gather}\label{a18}
\min_{\ell=1,2} \left( \frac{\partial \tR}{\partial t_\ell} +
D^{(\ell)}\tR(X_1,t_1,X_2,t_2)+g^{(\ell)}(X_1,t_1,X_2,t_2)\right)=0,
\end{gather}
where
\begin{gather*}
D^{(1)}\tR(X_1,t_1,X_2,t_2)=-\tR(X_1,t_1,X_2,t_2)X_1 t_1^{-1}\\
+\tR(X_1+1,t_1,X_2,t_2)(X_1+1) t_1^{-1},\\
D^{(2)}\tR(X_1,t_1,X_2,t_2)=-\tR(X_1,t_1,X_2,t_2)X_2 t_2^{-1}\\
+\tR(X_1,t_1,X_2+1,t_2)(X_2+1) t_2^{-1}.
\end{gather*}
Bayesian risk~\eqref{a3} is calculated by the formula
\begin{gather}\label{a19}
   R_T(\mu)=\tR(0,0,0,0).
\end{gather}

Note that partial differential equation at the same time describes
the evolution of  $\tR(X_1,t_1,X_2,t_2)$ and the strategy. The
strategy must choose $\ell$-th action if the $\ell$-th member in
the left-hand side of~\eqref{a18} has smaller value, in case of a
draw the choice of the action may be arbitrary.

%%%%%%%%%%%%%%%%%%%%%%%%%%%%%%%%%%%%%%%%%%

\end{document}